\documentclass{article}
\usepackage[british]{babel}
\usepackage{amsfonts}
\usepackage{amsbsy}
\usepackage{amssymb}

\usepackage{amsmath}

\DeclareMathOperator{\tr}{tr}
\DeclareMathOperator{\id}{id}
\DeclareMathOperator{\Aut}{Aut}

\usepackage{amsthm}
\theoremstyle{plain}
\newtheorem{theorem}{Theorem}
\theoremstyle{definition}
\newtheorem{definition}{Definition}

\begin{document}

\title{Spectral Expansions of Homogeneous and Isotropic Tensor-Valued Random Fields}

\author{Anatoliy Malyarenko\thanks{School of Education, Culture, and Communication, M\"{a}lardalen University, SE \textup{721 23} V\"{a}ster{\aa}s, Sweden, e-mail\texttt{anatoliy.malyarenko@mdh.se}} \and Martin Ostoja-Starzewski\thanks{Department of Mechanical Science and Engineering, University of Illinois at Urbana-Champaign, Urbana, IL, \textup{6181-2906}, USA, e-mail: \texttt{martinos@illinois.edu}}}

\maketitle

\begin{abstract}
We establish spectral expansions of homogeneous and isotropic random fields taking values in the $3$-dimensional Euclidean space $E^3$ and in the space $\mathsf{S}^2(E^3)$ of symmetric rank $2$ tensors over $E^3$. The former is a model of turbulent fluid velocity, while the latter is a model for the random stress tensor or the random conductivity tensor. We found a link between the theory of random fields and the theory of finite-dimensional convex compacta.
\end{abstract}

\section{Introduction}

Many random fields arising in continuum physics take values in linear spaces
of tensors over the \emph{space domain} $E^3$, the $3$-dimensional Euclidean space. For example, fluid velocity fields take values in the space of rank $1$ tensors. Stress, strain, rotation, and curvature-torsion fields take values in the space of rank $2$ tensors, while stiffness and compliance fields take values in the space of rank $4$ tensors. Their $n$-point correlation functions are shift-invariant. Under rotation, they transform according to an orthogonal representation of the orthogonal group $O(E^3)$.

To motivate the research in this direction, consider the differential form of Fourier's Law of thermal conduction which says that the local heat flux
density, $\mathbf{q}$ ($=q_i$), is equal to the product of thermal
conductivity, $k$, and the negative local temperature gradient, $-\boldsymbol{\nabla }T$ ($=T_{,i}$):%
\begin{equation}\label{eq:3}
\mathbf{q}=-k\boldsymbol{\nabla }T\textup{ \ \ \ or \ \ \ }q_{i}=-kT_{,i}
\end{equation}%
Here we use $\mathbf{A,B,...}$ for symbolic notation of a tensor, and $%
A_i,B_{ij}, \dots$\ for a subscript notation of tensors of 1st rank, 2nd rank, and so on...; a comma is used to indicate partial differentiation. Also, we use the (Einstein) summation convention (i.e. summing on the twice repeated subscript). The thermal conductivity, $k$, is often treated as a constant, though this is not always true. It may generally vary with temperature, which would make the heat conduction non-linear, and we do not consider it here.

In heat conduction the law of conservation of energy becomes%
\begin{equation}\label{eq:4}
\rho c\frac{\partial T}{\partial t}=-\nabla \cdot q_{i}\textup{ \ \ \ or \ \ \ }\rho c\frac{\partial T}{\partial t}=-q_{i},_{i}
\end{equation}%
where $\rho $\ is the mass density and $c$ is the specific heat capacity,
both assumed constant. Upon substitution of (\ref{eq:3}) into (\ref{eq:4}),
we find the heat conduction (or diffusion) equation%
\[
\frac{\partial T}{\partial t}=K\nabla ^{2}T\textup{ \ \ \ or \ \ \ }\rho c%
\frac{\partial T}{\partial t}=T,_{ii}
\]
where $K=k/\rho c$ and $\nabla ^{2}$\ is the Laplacian. In the case of
steady-state heat conduction, we get the Laplace equation%
\[
0=\nabla ^{2}T\textup{ \ \ \ or \ \ \ }0=T,_{ii}
\]

In nonuniform (i.e. inhomogeneous) media, $k$ varies with spatial location
over a spatial domain $D$, which is a subset of an $n$-dimensional Euclidean
space $E^{n}$ ($n=1,2$ or $3$). In general, we have an ensemble of
inhomogeneous media
\[
\{k(\omega ,\mathbf{x});\omega \in \Omega ,\mathbf{x}\in D\}
\]
so that $k(\omega ,\mathbf{x})$\ is a realisation of a random field (RF) $k$, $\Omega$ being the space of sample events. This is a good model if the
medium is piecewise constant (e.g. a polycrystal). However, for a
conductivity field to be random and have continuous realisations
\cite{Ostoja-Starzewski1999}, from microstructural considerations, it must be anisotropic at any given point $\mathbf{x}$. That is, the thermal conductivity $k$ must vary with orientation, and in this case $k$ is a second-rank tensor $k_{ij}$, and the Fourier law becomes
\[
\mathbf{q}=-\mathbf{k\cdot \nabla }T\textup{ \ \ \ or \ \ \ }q_{i}=-k_{ij}T_{,\textup{ }j}
\]
where $\mathbf{\cdot }$\ denotes a scalar product. The random medium is then
modelled by an ensemble of inhomogeneous, locally anisotropic media%
\[
\{k_{ij}(\omega ,\mathbf{x});\omega \in \Omega ,\mathbf{x}\in D\}
\]
such that, for any fixed $\omega $\ and $\mathbf{x}$, $k_{ij}$\ is a
positive definite, real-valued matrix. If we set $k_{ij}(\omega ,\mathbf{x}%
)=k(\omega ,\mathbf{x})$ $\delta _{ij}$, we recover a random medium with
locally isotropic realisations, but we note that in any random medium the
heat flux and temperature gradient are vector random fields. By virtue of the well known mathematical analogy, all the considerations above carry over to
in-plane states of stress, and, by extension to three dimensions, to stress
and strain fields as well as their connection via the 4th rank stiffness
tensor. Clearly, we need a more explicit way of representing and generating
vector and tensor random fields.

The statistical theory of isotropic turbulence was created by Sir Geoffrey Ingram Taylor \cite{Taylor} and developed further by numerous researchers. In particular, Robertson \cite{Robertson} proved that the correlation $R_{ij}$ between the $i$th component $u_i(\mathbf{x})$ of the velocity at $\mathbf{x}$ and the $j$th component $u_j(\mathbf{x}')$ at another point $\mathbf{x}'$ of the turbulent fluid is given by
\begin{equation}\label{eq:1}
R_{ij}=A\xi_i\xi_j+B\delta_{ij},
\end{equation}
where $\xi_i=x_i-x'_i$, and the coefficients $A$, $B$ are functions of the distance $\rho$ between $\mathbf{x}$ and $\mathbf{x}'$.

Lomakin \cite{Lomakin} considered the statistical theory of isotropic stress fields. He proved that the correlation $R_{ij\ell m}$ between the $ij$th component $\tau_{ij}(\mathbf{x})$ of the stress tensor at $\mathbf{x}$ and the $\ell m$th component $\tau_{\ell m}(\mathbf{x}')$ at another point $\mathbf{x}'$ of the body under deformation is
\begin{equation}\label{eq:2}
\begin{aligned}
R_{ij\ell m}&=a_1\delta_{ij}\delta_{\ell m}+a_2(\delta_{i\ell}\delta_{jm}
+\delta_{im}\delta_{j\ell})\\
&\quad+a_3(\xi_j\xi_{\ell}\delta_{im}
+\xi_i\xi_m\delta_{j\ell}+\xi_i\xi_{\ell}\delta_{jm}+\xi_j\xi_m\delta_{i\ell})\\
&\quad+a_4(\xi_i\xi_j\delta_{\ell m}+\xi_{\ell}\xi_m\delta_{ij})+a_5\xi_i\xi_j\xi_{\ell}\xi_m,
\end{aligned}
\end{equation}
where $a_1$, \dots, $a_5$ are functions of $\rho$.

In a different line of research, Yaglom \cite{Yaglom} proved that the correlation tensor (\ref{eq:1}) has the following spectral expansion:
\begin{equation}\label{eq:18}
\begin{aligned}
R_{ij}(\boldsymbol{\xi})&=\int^{\infty}_0
\left[\frac{j_1(\lambda\rho)}{\lambda\rho}\delta_{ij}-j_2(\lambda\rho)
\frac{\xi_i\xi_j}{\rho^2}\right]
\,\mathrm{d}\Phi_1(\lambda)\\
&\quad+\int^{\infty}_0\left[
\left(j_0(\lambda\rho)-\frac{j_1(\lambda\rho)}{\lambda\rho}\right)\delta_{ij}
+j_2(\lambda\rho)\frac{\xi_i\xi_j}{\rho^2}\right]
\,\mathrm{d}\Phi_2(\lambda),
\end{aligned}
\end{equation}
where $\Phi_1$ and $\Phi_2$ are two finite measures on $[0,\infty)$ with
\begin{equation}\label{eq:11}
\Phi_1(\{0\})=\Phi_2(\{0\})
\end{equation}
and where $j_i(t)$ are spherical Bessel functions. In particular, Robertson's functions $A(\rho)$ and $B(\rho)$ have the form
\[
\begin{aligned}
A(\rho)&=\frac{1}{\rho^2}\left(\int^{\infty}_0j_2(\lambda\rho)\,\mathrm{d}\Phi_2(\lambda)
-\int^{\infty}_0j_1(\lambda\rho)\,\mathrm{d}\Phi_1(\lambda)\right),\\
B(\rho)&=\int^{\infty}_0\frac{j_1(\lambda\rho)}{\lambda\rho}\,\mathrm{d}\Phi_1(\lambda)
+\int^{\infty}_0\left(j_0(\lambda\rho)-\frac{j_1(\lambda\rho)}{\lambda\rho}\right)
\,\mathrm{d}\Phi_2(\lambda).
\end{aligned}
\]

In this paper, we prove the spectral expansion of the correlation tensor (\ref{eq:2}) similar to that of Yaglom, and find the spectral expansions of both the turbulent fluid velocity field $u(\mathbf{x})$ and the stress field $\boldsymbol{\tau}(\mathbf{x})$ in terms of stochastic integrals with respect to orthogonal scattered random measures.

To achieve this goal, we first formulate our problem in mathematical language, introduce necessary notation and give the answer in Section~\ref{sec:2}. Then, we prove our results in Section~\ref{sec:3} and conclude in Section~\ref{sec:4}.

\section{Preliminaries}\label{sec:2}

Let $u(\mathbf{x})$ be the velocity of a turbulent fluid at a point $\mathbf{x}$ in the space domain $E^3$. Assume that $u(\mathbf{x})$ is a random field, i.e., a collection $\{\,u(\mathbf{x})\colon\mathbf{x}\in E^3\,\}$ of $E^3$-valued random vectors, defined on a probability space $(\Omega,\mathfrak{F},\mathsf{P})$. We suppose that the random field $u(\mathbf{x})$ is \emph{second-order}, i.e. $\mathsf{E}[\|u(\mathbf{x})\|^2]<\infty$, $\mathbf{x}\in E^3$, and \emph{mean-square continuous}, i.e., for any $\mathbf{x}_0\in E^3$ we have
\[
\lim_{\|\mathbf{x}-\mathbf{x}_0\|\to 0}\mathsf{E}[\|u(\mathbf{x})-u(\mathbf{x}_0)\|^2]=0.
\]
If one shifts the origin of the coordinate system by the vector $\mathbf{x}_0\in E^3$, the vector $u(\mathbf{x})$ does not change value. It follows that the random field $u(\mathbf{x})$ is \emph{wide-sense homogeneous}, i.e., its mean value $E(\mathbf{x}):=\mathsf{E}[u(\mathbf{x})]$ and correlation tensor $R(\mathbf{x},\mathbf{y}):=\mathsf{E}[(u(\mathbf{x})-E(\mathbf{x}))\otimes%
(u(\mathbf{y})-E(\mathbf{y}))]$
are shift-invariant: for any $\mathbf{x}_0\in E^3$ we have
\[
E(\mathbf{x}_0+\mathbf{x})=E(\mathbf{x}),\qquad R(\mathbf{x}_0+\mathbf{x},\mathbf{x}_0+\mathbf{y})=R(\mathbf{x},\mathbf{y}).
\]

Let $O(E^3)$ be the group of orthogonal linear transformations of the space $E^3$. Apply an arbitrary orthogonal transformation $k\in O(E^3)$ to the vector field $u(\mathbf{x})$. (Note that from now on $k$ denote an orthogonal transformation rather than the thermal conductivity tensor.) After the transformation $k$ the point $k^{-1}\mathbf{x}$ becomes the point $\mathbf{x}$. Evidently, the vector $u(k^{-1}\mathbf{x})$ is transformed by $k$ into $ku(k^{-1}\mathbf{x})$. It follows that for any positive integer $n$, for all distinct points $\mathbf{x}_1$, \dots, $\mathbf{x}_n\in E^3$, and for any $k\in O(E^3)$, the
random vectors $(u(\mathbf{x}_1),\dots,u(\mathbf{x}_n))^{\top}$ and $(ku(k^{-1}\mathbf{x}_1),\dots,ku(k^{-1}\mathbf{x}_n))^{\top}$ are identically distributed. Calculate the expectation of the transformed field:
\[
\mathsf{E}[ku(k^{-1}\mathbf{0})]=k\mathsf{E}[u(k^{-1}\mathbf{0})]=k\mathsf{E}
[u(\mathbf{0})].
\]
On the other hand,  $\mathsf{E}[ku(k^{-1}\mathbf{0})]=\mathsf{E}[u(\mathbf{0})]$. It follows that $k\mathsf{E}[u(\mathbf{0})]=\mathsf{E}[u(\mathbf{0})]$, $k\in O(E^3)$, therefore we have $E(\mathbf{x})=\mathbf{0}$.

Calculate the correlation function of the transformed field:
\[
\mathsf{E}[(ku(k^{-1}\mathbf{x}))\otimes(ku(k^{-1}\mathbf{y}))]=(k\otimes k) \mathsf{E}[u(k^{-1}\mathbf{x})\otimes u(k^{-1}\mathbf{y})].
\]
It follows that $R(k\boldsymbol{\xi})=(k\otimes k)R(\boldsymbol{\xi})$, where $\boldsymbol{\xi}=\mathbf{x}-\mathbf{y}$.

Let $\boldsymbol{\tau}(\mathbf{x})$ be the stress tensor of a deformable body. Assume that $\boldsymbol{\tau}(\mathbf{x})$ is a second-order mean-square continuous random field taking values in the space $\mathsf{S}^2(E^3)$ of symmetric rank $2$ tensors over $E^3$. Similar arguments prove that
\[
E(k\mathbf{x})=\mathsf{S}^2(k)E(\mathbf{x}),\qquad R(k\boldsymbol{\xi})=(\mathsf{S}^2(k)\otimes\mathsf{S}^2(k))R(\boldsymbol{\xi})
\]
for all $k\in O(E^3)$, where $\mathsf{S}^2(k)$ is the symmetric tensor square of the operator~$k$. Note that $k\mapsto\mathsf{S}^2(k)$ is an orthogonal representation of the group $O(E^3)$ in the space $L=\mathsf{S}^2(E^3)$.

We arrive at the following definition.  Let $r$ be a positive integer, and let $k\mapsto k^{\otimes r}$ be the orthogonal representation of the group $O(E^3)$ in the $r$th tensor power $(E^3)^{\otimes r}$ of the space $E^3$, let $L$ be an invariant subspace of the above representation, and let $U$ be the restriction of the above representation to $L$.

\begin{definition}
A random field $u(\mathbf{x})$, $\mathbf{x}\in E^3$ taking values in $L$ is called \emph{wide-sense isotropic} if
\begin{equation}\label{eq:6}
E(k\mathbf{x})=U(k)E(\mathbf{x}),\qquad R(k\mathbf{x},k\mathbf{y})=(U(k)\otimes U(k))R(\mathbf{x},\mathbf{y})
\end{equation}
for all $k\in O(E^3)$.
\end{definition}

In what follows, ``homogeneous random field'' always means ``wide-sense
homogeneous random field'', and ``isotropic random field'' always means ``wide-sense isotropic random field''.

In particular, in the case of the turbulent fluid velocity field we have $r=1$ and $L=E^3$, while in the case of the stress field we have $r=2$ and $L=\mathsf{S}^2(E^3)$. We would like to find the spectral expansion of both the correlation tensor of the stress field and the field itself, and to find the spectral expansion of the turbulent fluid velocity field.

Introduce the necessary notation. Let $\mathbb{K}$ be either the field $\mathbb{R}$ of real numbers or the field $\mathbb{C}$ of complex numbers. Let $V$ be a finite-dimensional vector space over $\mathbb{K}$, and let $\Aut V$ be the set of automorphisms of $V$. Let $K$ be a topological group with
identity element $e$. A \emph{representation} of the group $K$ in $V$ is a continuous homomorphism $U\colon K\to\Aut V$. A representation is called complex if $\mathbb{K}=\mathbb{C}$ and \emph{real} if $\mathbb{K}=\mathbb{R}$.

For example, let $K=SU(2)$ be the group of matrices of the following form:
\[
k=
\begin{pmatrix}
\alpha & \beta \\
-\overline{\beta} & \overline{\alpha}
\end{pmatrix}
,\qquad\alpha,\beta\in\mathbb{C},\quad|\alpha|^2+|\beta|^2=1.
\]
Let $V_{\ell}$ be the space of homogeneous polynomials of degree $2\ell$ in
two complex variables $\xi$ and $\eta$. The map
\[
U^{\ell}(k)f(\xi,\eta)=f(\overline{\alpha}\xi-\beta\eta,\overline{\beta}\xi
+\alpha\eta)
\]
is a complex representation of $K$.

Realise $E^3$ as the space of Hermitian matrices with zero trace in
$\mathbb{C}^2$. Such a matrix has the form
\[
A=
\begin{pmatrix}
x_0 & x_1+x_{-1}\mathrm{i} \\
x_1-x_{-1}\mathrm{i} & -x_0
\end{pmatrix}
,\qquad x_{-1},x_0,x_1\in\mathbb{R}.
\]
The map $A\mapsto k^{-1}Ak$, where $k$ is an element of the group $SU(2)$, is a rotation, i.e., an element of the group $K=SO(3)$ of orthogonal $3\times 3$ matrices with determinant $1$. The matrices $k$ and $-k$ determine the same rotation. Conversely, each rotation in $SO(3)$ corresponds to a pair of matrices $k$ and $-k$ in $SU(2)$.

If $\ell$ is a nonnegative integer, then the representation $U^{\ell}$
has the property $U^{\ell}(k)=U^{\ell}(-k)$. Therefore, $U^{\ell}$ is a complex representation of $SO(3)$.

Put
\[
(jf)(\xi,\eta)=\overline{f}(-\eta,\xi),\qquad f\in V_{\ell},
\]
where $\overline{f}$ is the polynomial with coefficients which are complex
conjugate to that of $f$. The map $j$ has the following properties
\[
j(c_1f_1+c_2f_2)=\overline{c_1}f_1+\overline{c_2}f_2,\qquad j^2=\id.
\]
In other words, $j$ is a \emph{real structure} on $V_{\ell}$.

Any complex vector space has many real structures. The structure $j$ has a
special property: it commutes with the representation $U^{\ell}$: $jU^{\ell}=U^{\ell}j$. We split $V_{\ell}$ into the subsets of eigenvectors with eigenvalues $+1$ and $-1$. These are vector spaces $V^+_{\ell}$ and $V^-_{\ell}$. The space $V^+_{\ell}$ is a vector space over $\mathbb{R}$. The restriction $U^{\ell,+}$ of the representation $U^{\ell}$ to $V^+_{\ell}$ is a real representation of $K=SO(3)$. Similarly for $V^-_{\ell}$.

Let $U_i$, $i=1$, $2$, be the representations of a topological group $K$ in the spaces $V_i$. A linear operator $A\colon V_1\to V_2$ is called an \emph{intertwining operator} if $AU_1=U_2A$. The representations $U_1$ and $U_2$ are called \emph{equivalent} if there exists an invertible intertwining operator $A\colon V_1\to V_2$.

For example, the multiplication by $\mathrm{i}$ is an invertible intertwining operator between equivalent real representations $U^{\ell,+}$ and
$U^{\ell,-}$. In what follows, we denote both representations by the same symbol $U^{\ell}$, and both spaces $V^+_{\ell}$ and $V^-_{\ell}$ by the same symbol $V_{\ell}$.

The \emph{direct sum} of representations $U_1$ and $U_2$ is the
representation $U_1\oplus U_2$ acting in the direct sum $V_1\oplus
V_2$ by
\[
(U_1\oplus U_2)(k)(\mathbf{x}\oplus\mathbf{y})=U_1(k)\mathbf{x}\oplus U_2(k)\mathbf{y}.
\]
Similarly, the \emph{tensor product} of representations $U_1$ and $U_2$ is the unique representation $U_1\otimes U_2$ acting on the elements of the form $\mathbf{x}\otimes\mathbf{y}$ of the tensor product $V_1\otimes V_2$ by
\[
(U_1\otimes U_2)(k)(\mathbf{x}\otimes\mathbf{y})=U_1(k)\mathbf{x}\otimes U_2(k)\mathbf{y}.
\]

If $K$ is a compact group, then it is possible to give $V$ an inner
product $(\boldsymbol{\cdot},\boldsymbol{\cdot})$ which is invariant under
$U$, i.e.,
\[
(U(k)\mathbf{x},U(k)\mathbf{y})=(\mathbf{x},\mathbf{y}),\qquad
k\in K,\quad\mathbf{x},\mathbf{y}\in V.
\]
Choose an orthonormal basis in $V$. Then, for a complex representation we can regard $U$ as taking values in the group $U(n)$ of unitary matrices of
order $n$, and we speak of a unitary representation. For a real
representation, $U$ takes values in the group $O(n)$ of orthogonal
matrices of order $n$, and we speak of an orthogonal representation.

A representation $U$ in a space $V$ is called \emph{reducible} if there
exists a proper (not equal to $\{\mathbf{0}\}$ or $V$) invariant
subspace $W$ with $U(k)(\mathbf{x})\in W$ for all $\mathbf{x}\in W$ and
all $k\in K$. Otherwise $U$ is called \emph{irreducible}.

For a compact group $K$, each representation is the direct sum of irreducible representations. Moreover, the above sum is unique in the following sense.
Let $V_i$ run over the spaces of all inequivalent irreducible representations of $K$, and let $m_i$, $n_i$ be nonnegative integers of which all but a
finite number are zero. Let $m_iV_i$ (resp. $n_iV_i$) be the direct sum of
$m_i$ (resp. $n_i$) copies of $V_i$. If $\oplus_im_iV_i$ is equivalent to
$\oplus_in_iV_i$, then $m_i=n_i$ for all $i$.

Above, we described the representatives of all equivalence classes of irreducible unitary representations of the group~$SU(2)$ and of all
irreducible unitary and orthogonal representations of the group~$SO(3)$. To describe the representatives of all equivalence classes of irreducible orthogonal representations of the group $O(3)$, we note that $O(3)$ is isomorphic to the direct product of the groups $SO(3)$ and $\mathbb{Z}_2$, and the subgroup $\mathbb{Z}_2$ is identical to $\{I,-I\}$, where $I$ is the identity matrix. Therefore, any irreducible orthogonal representation of the group $O(3)$ is isomorphic to the tensor product of irreducible orthogonal representations of the groups $SO(3)$ and $\mathbb{Z}_2$. In what follows, we denote by $U^{\ell,1}$ (resp. $U^{\ell,-1}$) the irreducible orthogonal representation of the group $O(3)$ with $U^{\ell,1}(-I)=\id$ (resp. $U^{\ell,-1}(-I)=-\id$), whose restriction to $SO(3)$ is equal to $U^{\ell}$. Note that the trivial representation of the group $O(3)$ is $U^{0,1}$, while the representation $k\mapsto k$ is $U^{1,-1}$.

For any representation $U$ of a compact group $K$ in a finite-dimensional space $V$, the spaces $\mathsf{S}^2(V)$ and $\wedge^2(V)$ of the symmetric and skew-symmetric rank $2$ tensors over $V$ are invariant
subspaces of the representation $U\otimes U$. Moreover, the
representation $U\otimes U$ is the direct sum of the corresponding
restrictions:
\[
U\otimes U=\mathsf{S}^2(U)\oplus\wedge^2(U).
\]

We introduce the basis $\mathbf{h}^{\ell}_m$, $-\ell\leqslant m\leqslant\ell$ in $V_{\ell}$ proposed by Gordienko \cite{Gordienko}. The Wigner $D$-functions, i.e., the matrix entries $D^{\ell,\pm 1}_{ij}(k)$ of the representations $U^{\ell,\pm 1}$ in the above basis are real-valued functions on the group $O(3)$.

The representation $U^{\ell_1,-1}\otimes U^{\ell_2,-1}$ is equivalent to the direct sum of the irreducible representations $U^{\ell,1}$, $|\ell_1-\ell_2|\leqslant\ell\leqslant\ell_1+\ell_2$. The transition from the \emph{uncoupled basis} $\{\,\mathbf{h}^{\ell}_m\colon|\ell_1-\ell_2|\leqslant\ell\leqslant\ell_1+\ell_2,-\ell\leqslant m\leqslant\ell\,\}$ to the \emph{coupled basis} $\{\,\mathbf{h}^{\ell_1}_{m_1}\otimes\mathbf{h}^{\ell_2}_{m_2}\colon-\ell_1\leqslant m_1\leqslant\ell_1,-\ell_2\leqslant m_2\leqslant\ell_2\,\}$ is performed by the \emph{Godunov--Gordienko coefficients}
\begin{equation}\label{eq:21}
\mathbf{h}^{\ell_1}_{m_1}\otimes\mathbf{h}^{\ell_2}_{m_2}
=\sum^{\ell_1+\ell_2}_{\ell=|\ell_1-\ell_2|}\sum^{\ell}_{m=-\ell}
g^{m[m_1,m_2]}_{\ell[\ell_1,\ell_2]}\mathbf{h}^{\ell}_m,
\end{equation}
introduced in \cite{Godunov}. By convention, we set $g^{m[m_1,m_2]}_{\ell[\ell_1,\ell_2]}=0$ if $\ell_1$, $\ell_2$, and $\ell$ do not satisfy the \emph{triangle condition} $|\ell_1-\ell_2|\leqslant\ell\leqslant\ell_1+\ell_2$.

Introduce the following notation:
\[
\begin{aligned}
b^{\ell'm'j}_{\ell mi,1}&=\mathrm{i}^{\ell-\ell'}\sqrt{(2\ell+1)(2\ell'+1)}
\left(\frac{1}{3}\delta_{ij}g^{0[m,m']}_{0[\ell,\ell']}g^{0[0,0]}_{0[\ell,\ell']}
-\frac{1}{5\sqrt{6}}g^{0[0,0]}_{2[\ell,\ell']}\sum^2_{n=-2}
g^{n[i,j]}_{2[1,1]}g^{-n[m,m']}_{2[\ell,\ell']}\right),\\
b^{\ell'm'j}_{\ell mi,2}&=\mathrm{i}^{\ell-\ell'}\sqrt{(2\ell+1)(2\ell'+1)}
\left(\frac{1}{3}\delta_{ij}g^{0[m,m']}_{0[\ell,\ell']}g^{0[0,0]}_{0[\ell,\ell']}
+\frac{\sqrt{2}}{5\sqrt{3}}g^{0[0,0]}_{2[\ell,\ell']}\sum^2_{n=-2}
g^{n[i,j]}_{2[1,1]}g^{-n[m,m']}_{2[\ell,\ell']}\right).
\end{aligned}
\]
Let $<$ be the lexicographic order on triples $(\ell,m,i)$, $\ell\geqslant  0$, $-\ell\leqslant m\leqslant\ell$, $-1\leqslant i\leqslant 1$. Let $L^1$ and $L^2$ be infinite lower triangular matrices from Cholesky factorisation of nonnegative-definite matrices $b^{\ell'm'j}_{\ell mi,1}$ and  $b^{\ell'm'j}_{\ell mi,2}$, constructed in \cite{Flinta}. Finally, let $Z^1_{\ell mi}$ and $Z^2_{\ell mi}$ be the set of centred uncorrelated random measures on $[0,\infty)$ with $\Phi_1$ being the control measure for $Z^1_{\ell mi}$ and $\Phi_2$ for $Z^2_{\ell mi}$.

The answers are given by the following theorems.

\begin{theorem}\label{th:1}
In the case of $V=\mathbb{R}^3$ and $U(k)=k$, the homogeneous and isotropic random field $u(\mathbf{x})$ has the form
\[
\begin{aligned}
u_i(r,\theta,\varphi)&=2\sqrt{\pi}\sum^{\infty}_{\ell=0}
\sum^{\ell}_{m=-\ell}\int^{\infty}_0j_{\ell}(\lambda r)\,
\mathrm{d}Z^{1'}_{\ell mi}(\lambda)S^m_{\ell}(\theta,\varphi)\\
&\quad+2\sqrt{\pi}\sum^{\infty}_{\ell=0}\sum^{\ell}_{m=-\ell}
\int^{\infty}_0j_{\ell}(\lambda r)\,\mathrm{d}Z^{2'}_{\ell mi}(\lambda)
S^m_{\ell}(\theta,\varphi),
\end{aligned}
\]
where
\[
Z^{k'}_{\ell mi}(A)=\sum_{(\ell',m',j)\leqslant(\ell,m,i)}L^k_{\ell mi,\ell'm'j}Z^k_{\ell'm'j}(A),
\]
with $k\in\{1,2\}$ and $A\in\mathfrak{B}([0,\infty))$.
\end{theorem}

\begin{table}
\caption{The functions $N_{nq}(\lambda,\rho)$}
\begin{tabular}{|l|l|l|}
\hline $n$ & $q$ & $N_{nq}(\lambda,\rho)$ \\
\hline
1 & 1 & $-\frac{2}{15}j_0(\lambda\rho)-\frac{4}{21}j_2(\lambda\rho)-\frac{2}{35}j_4(\lambda\rho)$ \\
1 & 2 & $\frac{1}{5}j_0(\lambda\rho)+\frac{1}{7}j_2(\lambda\rho)+\frac{26}{35}j_4(\lambda\rho)$ \\
1 & 3 & $-\frac{3}{14}j_2(\lambda\rho)+\frac{2}{7}j_4(\lambda\rho)$ \\
1 & 4 & $\frac{2}{7}(j_2(\lambda\rho)+j_4(\lambda\rho))$ \\
1 & 5 & $-2j_4(\lambda\rho)$ \\
2 & 1 & $-\frac{4}{45}j_0(\lambda\rho)+\frac{16}{63}j_2(\lambda\rho)+\frac{1}{105}j_4(\lambda\rho)$ \\
2 & 2 & $\frac{2}{15}j_0(\lambda\rho)-\frac{4}{21}j_2(\lambda\rho)-\frac{13}{105}j_4(\lambda\rho)$ \\
2 & 3 & $\frac{2}{7}j_2(\lambda\rho)-\frac{1}{21}j_4(\lambda\rho)$ \\
2 & 4 & $-\frac{8}{21}j_2(\lambda\rho)-\frac{1}{21}j_4(\lambda\rho)$ \\
2 & 5 & $\frac{1}{3}j_4(\lambda\rho)$ \\
3 & 1 & $\frac{2v_1(\lambda)+8v_2(\lambda)+1}{15}j_0(\lambda\rho)
+\frac{-8v_1(\lambda)+10v_2(\lambda)+2}{21}j_2(\lambda\rho)
+\frac{-v_1(\lambda)-4v_2(\lambda)+2}{70}j_4(\lambda\rho)$ \\
3 & 2 & $\frac{-v_1(\lambda)-4v_2(\lambda)+2}{30}j_0(\lambda\rho)
+\frac{-v_1(\lambda)-4v_2(\lambda)+2}{21}j_2(\lambda\rho)
+\frac{13(v_1(\lambda)+4v_2(\lambda)-2)}{70}j_4(\lambda\rho)$ \\
3 & 3 & $\frac{v_1(\lambda)+4v_2(\lambda)-2}{14}(j_2(\lambda\rho)
+j_4(\lambda\rho))$ \\
3 & 4 & $\frac{4v_1(\lambda)-5v_2(\lambda)-1}{7}j_2(\lambda\rho)
+\frac{v_1(\lambda)+4v_2(\lambda)-2}{14}j_4(\lambda\rho)$ \\
3 & 5 & $\frac{-v_1(\lambda)-4v_2(\lambda)+2}{2}j_4(\lambda\rho)$ \\
\hline
\end{tabular}
\label{tab:2}
\end{table}

Introduce the following notation.
\[
\begin{aligned}
L^1_{ij\ell m}(\boldsymbol{\xi})&=\delta_{ij}\delta_{\ell m},\\
L^2_{ij\ell
m}(\boldsymbol{\xi})&=\delta_{i\ell}\delta_{jm}+\delta_{im}\delta_{jl},\\
L^3_{ij\ell m}(\boldsymbol{\xi})&=\frac{\xi_j\xi_{\ell}}{\|\boldsymbol{\xi}\|^2}\delta_{im}
+\frac{\xi_i\xi_m}{\|\boldsymbol{\xi}\|^2}\delta_{j\ell}+\frac{\xi_i\xi_{\ell}}{\|
\boldsymbol{\xi}\|^2}\delta_{jm} +\frac{\xi_j\xi_m}{\|\boldsymbol{\xi}\|^2}\delta_{i\ell},\\
L^4_{ij\ell m}(\boldsymbol{\xi})&=\frac{\xi_i\xi_j}{\|\boldsymbol{\xi}\|^2}\delta_{\ell m}+\frac{\xi_{\ell}\xi_m}{\|\boldsymbol{\xi}\|^2}\delta_{ij},\\
L^5_{ij\ell
m}(\boldsymbol{\xi})&=\frac{\xi_i\xi_j\xi_{\ell}\xi_m}{\|\boldsymbol{\xi}\|^4}.
\end{aligned}
\]

\begin{theorem}\label{th:2}
In the case of $V=\mathsf{S}^2(\mathbb{R}^3)$ and $U(k)=\mathsf{S}^2(k)$, the expected value of the homogeneous and isotropic random field is
\[
E_{ij}(\mathbf{x})=C\delta_{ij},\qquad C\in\mathbb{R},
\]
while its correlation tensor has the spectral expansion
\begin{equation}\label{eq:31}
B_{ij\ell m}(\boldsymbol{\xi})=\sum^3_{n=1}\int^{\infty}_0\sum^5_{q=1}
N_{nq}(\lambda,\rho)L^q_{ij\ell m}(\boldsymbol{\xi})\,\mathrm{d}\Phi_n(\lambda),
\end{equation}
where the functions $N_{nq}(\lambda,\rho)$ are given in Table~\ref{tab:2}, $\Phi_n(\lambda)$ are three finite measures on $[0,\infty)$ with the following restriction: the atom $\Phi_3(\{0\})$ occupies at least $2/7$ of the sum of all three atoms, while the rest is divided between $\Phi_1(\{0\})$ and $\Phi_2(\{0\})$ in the proportion $1:\frac{3}{2}$. In Table~\ref{tab:2} $\mathbf{v}(\lambda)=(v_1(\lambda),v_2(\lambda))^{\top}$ is a $\Phi_3$-equivalence class of measurable functions taking values in the closed elliptic region $4(v_1(\lambda)-1/2)^2+8v^2_2(\lambda)\leqslant 1$.
\end{theorem}

In particular, we recover formula (\ref{eq:2}) with
\[
a_q(\rho)=\rho^{-s}\sum^3_{n=1}\int^{\infty}_0N_{nq}(\lambda,\rho)\,\mathrm{d}\Phi_n(\lambda),
\]
where $s=0$ for $q=1$, $2$, $s=2$ for $q=3$, $4$, and $s=4$ for $q=5$.

Introduce the following notation.
\[
\begin{aligned}
b^{u'w'\ell m}_{uwij,1}&=\mathrm{i}^{u-u'}\sqrt{(2u+1)(2u'+1)}\left(\frac{2}{5}
g^{0[w,w']}_{0[u,u']}g^{0[0,0]}_{0[u,u']}\sum^2_{n=-2}g^{n[i,j]}_{2[1,1]}
g^{n[\ell,m]}_{2[1,1]}\right.\\
&\quad+\frac{\sqrt{2}}{5\sqrt{7}}
g^{0[0,0]}_{2[u,u']}\sum^2_{n,q,t=-2}g^{-t[n,q]}_{2[2,2]}g^{n[i,j]}_{2[1,1]}
g^{q[\ell,m]}_{2[1,1]}g^{-t[w,w']}_{2[u,u']}\\
&\quad\left.-\frac{4\sqrt{2}}{3\sqrt{35}}
g^{0[0,0]}_{4[u,u']}\sum^2_{n,q=-2}\sum^4_{t=-4}g^{-t[n,q]}_{2[2,2]}
g^{n[i,j]}_{2[1,1]}g^{q[\ell,m]}_{2[1,1]}g^{-t[w,w']}_{4[u,u']}\right),\\
b^{u'w'\ell m}_{uwij,2}&=\mathrm{i}^{u-u'}\sqrt{(2u+1)(2u'+1)}\left(\frac{4}{15}
g^{0[w,w']}_{0[u,u']}g^{0[0,0]}_{0[u,u']}\sum^2_{n=-2}g^{n[i,j]}_{2[1,1]}
g^{n[\ell,m]}_{2[1,1]}\right.\\
&\quad-\frac{4\sqrt{2}}{15\sqrt{7}}
g^{0[0,0]}_{2[u,u']}\sum^2_{n,q,t=-2}g^{-t[n,q]}_{2[2,2]}g^{n[i,j]}_{2[1,1]}
g^{q[\ell,m]}_{2[1,1]}g^{-t[w,w']}_{2[u,u']}\\
&\quad\left.+\frac{2\sqrt{2}}{27\sqrt{35}}
g^{0[0,0]}_{4[u,u']}\sum^2_{n,q=-2}\sum^4_{t=-4}g^{-t[n,q]}_{2[2,2]}
g^{n[i,j]}_{2[1,1]}g^{q[\ell,m]}_{2[1,1]}g^{-t[w,w']}_{4[u,u']}\right),\\
b^{u'w'\ell m}_{uwij,3}(\lambda)&=\mathrm{i}^{u-u'}\sqrt{(2u+1)(2u'+1)}\left(
\frac{v_1(\lambda)+4v_2(\lambda)+1}{9}\delta_{ij}\delta_{\ell m}g^{0[w,w']}_{0[u,u']}g^{0[0,0]}_{0[u,u']}\right.\\
&\quad+\frac{-v_1(\lambda)-4v_2(\lambda)+2}{15}
g^{0[w,w']}_{0[u,u']}g^{0[0,0]}_{0[u,u']}\sum^2_{n=-2}g^{n[i,j]}_{2[1,1]}
g^{n[\ell,m]}_{2[1,1]}\\
&\quad+\frac{-4v_1(\lambda)+2v_2(\lambda)+2}{15\sqrt{6}}g^{0[0,0]}_{2[u,u']}
\sum^2_{t=-2}g^{-t[w,w']}_{2[u,u']}(\delta_{ij}g^{t[\ell,m]}_{2[1,1]}+\delta_{\ell m}g^{t[i,j]}_{2[1,1]})\\
&\quad+\frac{\sqrt{2}[-v_1(\lambda)-4v_2(\lambda)+2]}{15\sqrt{7}}
g^{0[0,0]}_{2[u,u']}\sum^2_{n,q,t=-2}g^{-t[n,q]}_{2[2,2]}g^{n[i,j]}_{2[1,1]}
g^{q[\ell,m]}_{2[1,1]}g^{-t[w,w']}_{2[u,u']}\\
&\quad\left.+\frac{\sqrt{2}[-v_1(\lambda)-4v_2(\lambda)+2]}{9\sqrt{35}}
g^{0[0,0]}_{4[u,u']}\sum^2_{n,q=-2}\sum^4_{t=-4}g^{-t[n,q]}_{2[2,2]}
g^{n[i,j]}_{2[1,1]}g^{q[\ell,m]}_{2[1,1]}g^{-t[w,w']}_{4[u,u']}\right).
\end{aligned}
\]
Let $<$ be the lexicographic order on quadruples $(u,w,i,j)$, $u\geqslant 0$, $-u\leqslant w\leqslant u$, $-1\leqslant i\leqslant 1$, $-1\leqslant j\leqslant 1$. Let $L^1$, $L^2$ and $L^3(\lambda)$ be infinite
lower triangular matrices from Cholesky factorisation of nonnegative-definite matrices $b^{u'w'\ell m}_{uwij,1}$, $b^{u'w'\ell m}_{uwij,2}$ and  $b^{u'w'\ell m}_{uwij,3}(\lambda)$, constructed in \cite{Flinta}. Finally, let $Z^1_{uwij}$, $Z^2_{uwij}$, and $Z^3_{uwij}$ be the set of centred uncorrelated random measures on $[0,\infty)$ with $\Phi_n$ being the control measure for $Z^n_{uwij}$, $1\leqslant n\leqslant 3$.

\begin{theorem}\label{th:3}
In the case of $V=\mathsf{S}^2(\mathbb{R}^3)$ and $\theta(k)=\mathsf{S}^2(k)$, the homogeneous and isotropic random field $\boldsymbol{\tau}(\mathbf{x})$ has the form
\[
\begin{aligned}
\tau_{ij}(r,\theta,\varphi)&=C\delta_{ij}+2\sqrt{\pi}\sum^{\infty}_{u=0}\sum^u_{w=-u}
\int^{\infty}_0j_u(\lambda r)\,\mathrm{d}Z^{1'}_{uwij}(\lambda)S^w_u(\theta,\varphi)\\
&\quad+2\sqrt{\pi}\sum^{\infty}_{u=0}\sum^u_{w=-u}
\int^{\infty}_0j_u(\lambda r)\,\mathrm{d}Z^{2'}_{uwij}(\lambda)S^w_u(\theta,\varphi)\\
&\quad+2\sqrt{\pi}\sum^{\infty}_{u=0}\sum^u_{w=-u}
\int^{\infty}_0j_u(\lambda r)\sum_{(u',w',i',j')\leqslant(u,w,i,j)}L^3_{uwij,u'w'i'j'}(\lambda)\\
&\quad\times\mathrm{d}Z^3_{uwij}(\lambda)S^w_u(\theta,\varphi),
\end{aligned}
\]
where
\begin{equation}\label{eq:38}
Z^{n'}_{uwij}(A)=\sum_{(u',w',i',j')\leqslant(u,w,i,j)}L^n_{uwij,u'w'i'j'}
Z^k_{u'w'i'j'}(A),
\end{equation}
with $1\leqslant k\leqslant 3$ and $A\in\mathfrak{B}([0,\infty))$.
\end{theorem}

\section{Proofs}\label{sec:3}

Proofs of Theorems~\ref{th:1}--\ref{th:3} have a common part that is applicable to a general homogeneous and isotropic random field $u(\mathbf{x})$.

Let the representation $U$ be the direct sum of $\ell_0$ copies of the irreducible orthogonal representation $U^{0,(-1)^r}$, $\ell_1$ copies of the representation $U^{1,(-1)^r}$, \dots, $\ell_r$ copies of the representation $U^{r,(-1)^r}$. Let $T^{i,j,n}_{m_1\cdots m_r}$, $-i\leqslant n\leqslant i$, be the vectors of the Gordienko basis of the space where the $j$th copy of the representation $U^{i,(-1)^r}$ acts. The rank $r$ tensors
\[
T^{i,j,n}_{m_1\cdots m_r}\colon 0\leqslant i\leqslant r,1\leqslant j\leqslant \ell_i,-i\leqslant n\leqslant i
\]
constitute the uncoupled basis of the space $L$. In the first equation in (\ref{eq:6}), put $\mathbf{x}=\mathbf{0}$. We obtain $E(\mathbf{0})=U(k)E(\mathbf{0})$, $k\in O(3)$. In other words, $E(\mathbf{x})$ lies in the space where the direct sum of $\ell_0$ copies of the trivial representation $U^{0,1}$ acts. This space may have positive dimension if $r$ is even and $\ell_0>0$. In this case we obtain
\begin{equation}\label{eq:22}
E(\mathbf{x})=\sum^{\ell_0}_{j=1}C_jT^{0,j,0}_{m_1\cdots m_r},\qquad C_j\in\mathbb{R}.
\end{equation}

Let $W=L\oplus\mathrm{i}L$ be the complexification of the space $L$. It is known (cf.\cite[Theorem~2 and Remark~1]{Yaglom}) that equation
\begin{equation}\label{eq:5}
R(\boldsymbol{\xi})=\int_{\hat{\mathbb{R}}^3}
\mathrm{e}^{\mathrm{i}(\mathbf{p},\boldsymbol{\xi})}\,\mathrm{d}F(\mathbf{p}),
\end{equation}
where $\hat{\mathbb{R}}^3$ is the \emph{wavenumber domain}, establishes a one-to-one correspondence between correlation tensors $R(\boldsymbol{\xi})$ of homogeneous $W$-valued random fields and measures $F$ defined on the Borel $\sigma$-field $\mathfrak{B}(\hat{\mathbb{R}}^3)$ and taking values in the set of Hermitian nonnegative-definite linear operators in $W$. The set of all Hermitian operators in $L\oplus\mathrm{i}L$ is $\mathsf{S}^2(L)\oplus\mathrm{i}\wedge^2(L)$. Let $J$ be the linear operator in the above space acting by
\[
J(S+\mathrm{i}\wedge)=S-\mathrm{i}\wedge.
\]
If the random field takes values in $L$, then for any $A\in\mathfrak{B}(\hat{\mathbb{R}}^3)$ we have
\begin{equation}\label{eq:9}
F(-A)=JF(A),
\end{equation}
where $-A=\{\,-\mathbf{p}\colon\mathbf{p}\in A\,\}$.

Let $\sigma$ be the following measure:
\[
\sigma(A)=\tr[F(A)],\qquad A\in\mathfrak{B}(\hat{\mathbb{R}}^3),
\]
where $\tr$ denote the trace of a matrix. By \cite{Berezanskii}, the measure $F$ is absolutely continuous with respect to $\sigma$, and the density
$f(\mathbf{p})=dF(\mathbf{p})/d\sigma(\mathbf{p})$ is a
measurable function on $\hat{\mathbb{R}}^3$ taking values in the set of
Hermitian nonnegative-definite operators in the space $W$ with unit trace. Thus, equation (\ref{eq:5}) may be written as
\begin{equation}\label{eq:7}
R(\boldsymbol{\xi})=\int_{\hat{\mathbb{R}}^3}
e^{\mathrm{i}(\mathbf{p},\boldsymbol{\xi})}f(\mathbf{p})\,\mathrm{d}\sigma(\mathbf{p}).
\end{equation}

We calculate the expression $R(k\boldsymbol{\xi})$ by two different methods. On the one hand, by the second equation in (\ref{eq:6}),
\[
\begin{aligned}
R(k\boldsymbol{\xi})&=(U(k)\otimes U(k))R(\boldsymbol{\xi})\\
&=(U(k)\otimes U(k))\int_{\hat{\mathbb{R}}^3}e^{\mathrm{i}(\mathbf{p},
\boldsymbol{\xi})}f(\mathbf{p})\,\mathrm{d}\sigma(\mathbf{p})\\
&=\int_{\hat{\mathbb{R}}^3}e^{\mathrm{i}(\mathbf{p},\boldsymbol{\xi})}
(U(k)\otimes U(k))f(\mathbf{p})\,\mathrm{d}\sigma(\mathbf{p}),
\end{aligned}
\]
because integration commutes with continuous linear operators. On the other
hand, we have
$e^{\mathrm{i}(\mathbf{p},k\boldsymbol{\xi})}=e^{\mathrm{i}(k^{-1}\mathbf{p},\boldsymbol{\xi})}$.
Then, by (\ref{eq:7}),
\[
\begin{aligned}
R(k\boldsymbol{\xi})&=\int_{\hat{\mathbb{R}}^3}
e^{\mathrm{i}(\mathbf{p},k\boldsymbol{\xi})}f(\mathbf{p})\,\mathrm{d}\sigma(\mathbf{p})\\
&=&\int_{\hat{\mathbb{R}}^3}e^{\mathrm{i}(k^{-1}\mathbf{p},\boldsymbol{\xi})}
f(\mathbf{p})\,\mathrm{d}\sigma(\mathbf{p})\\
&=\int_{\hat{\mathbb{R}}^3}e^{\mathrm{i}(\mathbf{p},\boldsymbol{\xi})}
f(k\mathbf{p})\,\mathrm{d}\sigma(k\mathbf{p}).
\end{aligned}
\]
In the last display we denote $k^{-1}\mathbf{p}$ again by $\mathbf{p}$.
Because the expansion (\ref{eq:7}) is unique, we have, for each $k\in O(3)$
and for each $A\in\mathfrak{B}(\hat{\mathbb{R}}^3)$,
\begin{equation}\label{eq:8}
f(k\mathbf{p})=(U(k)\otimes U(k))f(\mathbf{p}),\qquad\sigma(kA)
=\sigma(A).
\end{equation}

Let $\mathrm{d}\Omega$ be the Lebesgue measure on the unit sphere $S^2\subset\hat{\mathbb{R}}^3$. The measure $\sigma$ satisfying the second part of (\ref{eq:8}), has the form
\begin{equation}\label{eq:12}
\mathrm{d}\sigma=(4\pi)^{-1}\mathrm{d}\Omega\,\mathrm{d}\mu(\lambda),
\end{equation}
where $\mu$ is a finite measure on $[0,\infty)$.

In the first equation of (\ref{eq:8}), put $k=-I$. We obtain $f(-\mathbf{p}=f(\mathbf{p})$. It follow from (\ref{eq:8}) and (\ref{eq:9}) that $f$ takes values in the subspace $\mathsf{S}^2(L)$. The first equation in (\ref{eq:8}) takes the form
\begin{equation}\label{eq:14}
f(k\mathbf{p})=\mathsf{S}^2(U(k))f(\mathbf{p}).
\end{equation}

Let the representation $\mathsf{S}^2(U)$ be be the direct sum of $\ell'_0$ copies of the irreducible orthogonal representation $U^{0,1}$, $\ell'_1$ copies of the representation $U^{1,1}$, \dots, $\ell'_{2r}$ copies of the representation $U^{2r,1}$. Let $T^{i,j,n}_{m_1\cdots m_{2r}}$, $-i\leqslant n\leqslant i$, be the vectors of the Gordienko basis of the space where the $j$th copy of the representation $U^{i,1}$ acts. The rank $r$ tensors
\[
T^{i,j,n}_{m_1\cdots m_{2r}}\colon 0\leqslant i\leqslant 2r,1\leqslant j\leqslant \ell'_i,-i\leqslant n\leqslant i
\]
constitute the uncoupled basis of the space $\mathsf{S}^2(L)$.

Let $(\lambda,\theta_{\mathbf{p}},\varphi_{\mathbf{p}})$ be the spherical coordinates in the wavenumber domain. Let $f^{i,j,n}(\lambda)$ be the value of the linear form $f(\lambda,0,0)$ on the tensor $T^{i,j,n}_{m_1\cdots m_{2r}}$. Let $\mathbf{f}^{i,j}(\lambda)\in\mathbb{R}^{2i+1}$ be the vector with coordinates $f^{i,j,n}(\lambda)$, $-i\leqslant n\leqslant i$. The stationary subgroup of the point $(0,0,0)$ is $O(3)$. It follows from (\ref{eq:14}) that
\[
\mathbf{f}^{i,j}(0)=U^{i,1}(k)\mathbf{f}^{i,j}(0),\qquad k\in O(3).
\]
It follows that $\mathbf{f}^{i,j}(0)=\mathbf{0}$, if $i\geqslant 1$.

For $\lambda>0$, the stationary subgroup of the point $(\lambda,0,0)$ is $O(2)$. By \cite[Claim~8.3]{Klimyk}, the restriction of the representation $U^{i,1}$ to the group $O(2)$ contains the trivial representation of $O(2)$ if and only if $i$ is even. It follows that only the functions $f^{2i,j,0}(\lambda)$ may be nonzero. By linearity, the matrix entries of the matrix $f(\lambda)$, i.e., the values of the linear functional $f(\lambda)$ on the tensor $\mathbf{e}_{m_1}\otimes\mathbf{e}_{m_2}\otimes\cdots\otimes\mathbf{e}_{m_{2r}}$, are as follows
\begin{equation}\label{eq:15}
f_{m_1\cdots m_{2r}}(\lambda)=\sum^r_{i=0}\sum^{\ell'_{2i}}_{j=1}
T^{2i,j,0}_{m_1\cdots m_{2r}}f^{2i,j,0}(\lambda).
\end{equation}
In other words, for all $\lambda\geqslant 0$, the matrix $f(\lambda)$ lies in the intersection $\mathcal{C}$ of the convex compact set of all nonnegative-definite matrices with unit trace and convex linear subspaces of the space of all matrices. It follows that $\mathcal{C}$ is a convex compact set. The structure of the extreme points of $\mathcal{C}$ will be analysed for each of the cases separately.

The value of the linear functional $f(\lambda,\theta_{\mathbf{p}},\varphi_{\mathbf{p}})$ on the tensor $T^{i,j,n}_{m_1\cdots m_{2r}}$ is calculated by
\[
f^{i,j,n}(\lambda,\theta_{\mathbf{p}},\varphi_{\mathbf{p}})=
D^{2i}_{n0}(\theta_{\mathbf{p}},\varphi_{\mathbf{p}})f^{i,j,0}(\lambda),
\]
which follows from (\ref{eq:14}). Again by linearity we obtain
\begin{equation}\label{eq:17}
f_{m_1\cdots m_{2r}}(\lambda,\theta_{\mathbf{p}},\varphi_{\mathbf{p}})=
\sum^r_{i=0}\sum^{\ell'_{2i}}_{j=1}\sum^{2i}_{n=-2i}T^{2i,j,n}_{m_1\cdots m_{2r}}D^{2i}_{n0}(\theta_{\mathbf{p}},\varphi_{\mathbf{p}})f^{2i,j,0}(\lambda).
\end{equation}
Finally, we introduce the following notation
\[
M^{2i,j}_{m_1\cdots m_{2r}}(\mathbf{p})=\sum^{2i}_{n=-2i}T^{2i,j,n}_{m_1\cdots m_{2r}}D^{2i}_{n0}(\theta_{\mathbf{p}},\varphi_{\mathbf{p}}).
\]

\subsection{Proof of Theorem~\ref{th:1}}

In Theorem~\ref{th:1} we have $r=1$ and $L=E^3$. The expected value of the random field is equal to $\mathbf{0}$ because $r$ is odd. The uncoupled basis in $\mathsf{S}^2(E^3)$ is constituted by the tensors $T^{0,1,0}_{ij}=g^{0[i,j]}_{0[1,1]}$ and $T^{2,1,n}_{ij}=g^{n[i,j]}_{2[1,1]}$. Using (\ref{eq:15}) and the values of the Godunov--Gordienko coefficients calculated in \cite{Godunov,Malyarenko}, we obtain
\[
f_{ij}(\lambda)=\left[\frac{2}{\sqrt{3}}f^{0,1,0}(\lambda)
-\frac{\sqrt{2}}{\sqrt{3}}f^{2,1,0}(\lambda)\right]D^1
+\left[\frac{1}{\sqrt{3}}f^{0,1,0}(\lambda)
+\frac{\sqrt{2}}{\sqrt{3}}f^{2,1,0}(\lambda)\right]D^2,
\]
where $D^1$ is the $3\times 3$ matrix with nonzero entries $D^1_{-1-1}=D^1_{11}=1/2$, while $D^2$ has the only nonzero entry $D^2_{00}=1$. In other words, the set $\mathcal{C}$ is the interval with extreme points $D^1$ and $D^2$, while
\[
u_1(\lambda)=\frac{2}{\sqrt{3}}f^{0,1,0}(\lambda)-\frac{\sqrt{2}}{\sqrt{3}}f^{2,1,0}(\lambda),
\qquad u_2(\lambda)=\frac{1}{\sqrt{3}}f^{0,1,0}(\lambda)
+\frac{\sqrt{2}}{\sqrt{3}}f^{2,1,0}(\lambda)
\]
are affine coordinates in the one-dimensional simplex $C$ with $u_1(\lambda)\geqslant 0$, $u_2(\lambda)\geqslant 0$, and $u_1(\lambda)+u_2(\lambda)=1$. Moreover, we have $f^{2,1,0}(0)=0$ and
\begin{equation}\label{eq:10}
u_1(0)=\frac{2}{\sqrt{3}}f^{0,1,0}(0),\qquad u_2(0)=\frac{1}{\sqrt{3}}f^{0,1,0}(0).
\end{equation}

The functions $f^{i,j,0}(\lambda)$ are expressed in terms of $u_1(\lambda)$ and $u_2(\lambda)$ as follows:
\[
f^{0,1,0}(\lambda)=\frac{1}{\sqrt{3}}u_1(\lambda)+\frac{1}{\sqrt{3}}u_2,
\qquad f^{2,1,0}(\lambda)=-\frac{1}{\sqrt{6}}u_1(\lambda)+\frac{\sqrt{2}}{\sqrt{3}}u_2(\lambda).
\]
Substitute these values to (\ref{eq:17}). We obtain
\[
f_{ij}(\mathbf{p})=\left[\frac{1}{\sqrt{3}}M^{0,1}_{ij}(\mathbf{p})
-\frac{1}{\sqrt{6}}M^{2,1}_{ij}(\mathbf{p})\right]u_1(\lambda)
+\left[\frac{1}{\sqrt{3}}M^{0,1}_{ij}(\mathbf{p})
+\frac{\sqrt{2}}{\sqrt{3}}M^{2,1}_{ij}(\mathbf{p})\right]u_2(\lambda).
\]
or
\[
\begin{aligned}
f_{ij}(\mathbf{p})&=\left(\frac{1}{3}
\delta_{ij}D^0_{00}(\theta_{\mathbf{p}},\varphi_{\mathbf{p}})
-\frac{1}{\sqrt{6}}\sum^2_{m=-2}g^{m[i,j]}_{2[1,1]}D^{2}_{m0}
(\theta_{\mathbf{p}},\varphi_{\mathbf{p}})\right)u_1(\lambda)\\
&\quad+\left(\frac{1}{3}\delta_{ij}D^0_{00}(\theta_{\mathbf{p}},\varphi_{\mathbf{p}})
+\frac{\sqrt{2}}{\sqrt{3}}\sum^2_{m=-2}g^{m[i,j]}_{2[1,1]}
D^2_{m0}(\theta_{\mathbf{p}},\varphi_{\mathbf{p}})\right)u_2(\lambda).
\end{aligned}
\]

Substitute this formula and (\ref{eq:12}) to (\ref{eq:7}). We obtain
\[
\begin{aligned}
R_{ij}(\boldsymbol{\xi})&=\frac{1}{4\pi}\int_{\hat{\mathbb{R}}^3}
e^{\mathrm{i}(\mathbf{p},\boldsymbol{\xi})}\left(\frac{1}{3}
\delta_{ij}D^0_{00}(\theta_{\mathbf{p}},\varphi_{\mathbf{p}})\right.\\
&\quad\left.-\frac{1}{\sqrt{6}}\sum^2_{m=-2}g^{m[i,j]}_{2[1,1]}D^2_{m0}
(\theta_{\mathbf{p}},\varphi_{\mathbf{p}})\right)u_1(\lambda)\,
\mathrm{d}\Omega\,\mathrm{d}\mu(\lambda)\\
&\quad+\frac{1}{4\pi}\int_{\hat{\mathbb{R}}^3}
e^{\mathrm{i}(\mathbf{p},\boldsymbol{\xi})}\left(\frac{1}{3}\delta_{ij}
D^0_{00}(\theta_{\mathbf{p}},\varphi_{\mathbf{p}})\right.\\
&\quad\left.+\frac{\sqrt{2}}{\sqrt{3}}\sum^2_{m=-2}g^{m[i,j]}_{2[1,1]}
D^2_{m0}(\theta_{\mathbf{p}},\varphi_{\mathbf{p}})\right)u_2(\lambda)\,
\mathrm{d}\Omega\,\mathrm{d}\mu(\lambda).
\end{aligned}
\]

The norm of the function $D^{\ell}_{m0}(\theta_{\mathbf{p}},\varphi_{\mathbf{p}})$ in the space $L^2(S^2,\mathrm{d}\Omega)$ is not equal to $1$, while the norm of the \emph{real-valued spherical harmonic}
\begin{equation}\label{eq:16}
S^m_{\ell}(\theta_{\mathbf{p}},\varphi_{\mathbf{p}})=\sqrt{\frac{2\ell+1}{4\pi}}
D^{\ell}_{-m0}(\theta_{\mathbf{p}},\varphi_{\mathbf{p}})
\end{equation}
is equal to $1$. In terms of spherical harmonics, we have
\begin{equation}\label{eq:19}
\begin{aligned}
R_{ij}(\boldsymbol{\xi})&=\frac{1}{\sqrt{4\pi}}\int_{\hat{\mathbb{R}}^3}
e^{\mathrm{i}(\mathbf{p},\boldsymbol{\xi})}\left(\frac{1}{3}
\delta_{ij}S^0_0(\theta_{\mathbf{p}},\varphi_{\mathbf{p}})\right.\\
&\quad\left.-\frac{1}{\sqrt{30}}\sum^2_{m=-2}g^{m[i,j]}_{2[1,1]}S^{-m}_2
(\theta_{\mathbf{p}},\varphi_{\mathbf{p}})\right)\,
\mathrm{d}\Omega\,\mathrm{d}\Phi_1(\lambda)\\
&\quad+\frac{1}{\sqrt{4\pi}}\int_{\hat{\mathbb{R}}^3}
e^{\mathrm{i}(\mathbf{p},\boldsymbol{\xi})}\left(\frac{1}{3}\delta_{ij}
S^0_0(\theta_{\mathbf{p}},\varphi_{\mathbf{p}})\right.\\
&\quad\left.+\frac{\sqrt{2}}{\sqrt{15}}\sum^2_{m=-2}g^{m[i,j]}_{2[1,1]}
S^{-m}_2(\theta_{\mathbf{p}},\varphi_{\mathbf{p}})\right)\,
\mathrm{d}\Omega\,\mathrm{d}\Phi_2(\lambda),
\end{aligned}
\end{equation}
where we introduced notation $\mathrm{d}\Phi_j(\lambda)=u_j(\lambda)\,\mathrm{d}\mu(\lambda)$, $j=1$, $2$. It follows from (\ref{eq:10}) that $\Phi_1(\{0\})=\frac{2}{3}\mu(\{0\})$ and $\Phi_2(\{0\})=\frac{1}{3}\mu(\{0\})$. In other words, $\Phi_1(\{0\})=2\Phi_2(\{0\})$, which differs from (\ref{eq:11}).

This may be explained as follows. The values of $M^{0,1}(\mathbf{p})$ and $M^{2,1}_{ij}(\mathbf{p})$ were calculated in \cite{Malyarenko} as
\begin{equation}\label{eq:13}
M^{0,1}_{ij}(\mathbf{p})=\frac{1}{\sqrt{3}}\delta_{ij},\qquad M^{2,1}_{ij}(\mathbf{p})=\frac{\sqrt{3}}{\sqrt{2}}\frac{p_ip_j}{\|\mathbf{p}\|^2}
-\frac{1}{\sqrt{6}}\delta_{ij}.
\end{equation}
Therefore we have
\[
f_{ij}(\mathbf{p})=\frac{1}{2}\left(\delta_{ij}-\frac{p_ip_j}{\|\mathbf{p}\|^2}\right)
u_1(\lambda)+\frac{p_ip_j}{\|\mathbf{p}\|^2}u_2(\lambda).
\]
In \cite{Yaglom} Yaglom uses invariant theory to prove the formula
\[
f_{ij}(\mathbf{p})=\left(\delta_{ij}-\frac{p_ip_j}{\|\mathbf{p}\|^2}\right)
u_1(\lambda)+\frac{p_ip_j}{\|\mathbf{p}\|^2}u_2(\lambda),
\]
i.e., his function $u_1(\lambda)$ is twice less than our one, hence the difference.

To calculate the inner integral in (\ref{eq:19}), use the following \emph{plane wave expansion}:
\begin{equation}\label{eq:20}
e^{\mathrm{i}(\mathbf{p},\boldsymbol{\xi})}=4\pi\sum^{\infty}_{\ell=0}
\mathrm{i}^{\ell}j_{\ell}(\lambda\rho)\sum^{\ell}_{m=-\ell}S^m_{\ell}
(\theta_{\boldsymbol{\xi}},\varphi_{\boldsymbol{\xi}})
S^m_{\ell}(\theta_{\mathbf{p}},\varphi_{\mathbf{p}}),
\end{equation}
where $(\rho,\theta_{\boldsymbol{\xi}},\varphi_{\boldsymbol{\xi}})$ are the spherical coordinates in the space domain. The spectral expansion takes the form
\[
\begin{aligned}
R_{ij}(\boldsymbol{\xi})&=\sqrt{4\pi}\int^{\infty}_0
\left(\frac{1}{3}\delta_{ij}j_0(\lambda\rho)
S^0_0(\theta_{\boldsymbol{\xi}},\varphi_{\boldsymbol{\xi}})\right.\\
&\quad\left.+\frac{1}{\sqrt{30}}j_2(\lambda\rho)\sum^2_{m=-2}g^{m[i,j]}_{2[1,1]}S^{-m}_2
(\theta_{\boldsymbol{\xi}},\varphi_{\boldsymbol{\xi}})\right)
\,\mathrm{d}\Phi_1(\lambda)\\
&\quad+\sqrt{4\pi}\int^{\infty}_0\left(\frac{1}{3}\delta_{ij}j_0(\lambda\rho)
S^0_0(\theta_{\boldsymbol{\xi}},\varphi_{\boldsymbol{\xi}})\right.\\
&\quad\left.-\frac{\sqrt{2}}{\sqrt{15}}j_2(\lambda\rho)\sum^2_{m=-2}g^{m[i,j]}_{2[1,1]}
S^{-m}_2(\theta_{\boldsymbol{\xi}},\varphi_{\boldsymbol{\xi}})\right)\,
\mathrm{d}\Phi_2(\lambda).
\end{aligned}
\]
Using (\ref{eq:16}) and (\ref{eq:13}), we obtain
\[
\begin{aligned}
R_{ij}(\boldsymbol{\xi})&=\int^{\infty}_0
\left[\left(\frac{1}{3}j_0(\lambda\rho)-\frac{1}{6}j_2(\lambda\rho)\right)
\delta_{ij}+\frac{1}{2}j_2(\lambda\rho)\frac{\xi_i\xi_j}{\|\boldsymbol{\xi}\|^2}\right]
\,\mathrm{d}\Phi_1(\lambda)\nonumber\\
&\quad+\int^{\infty}_0\left[
\left(\frac{1}{3}j_0(\lambda\rho)+\frac{1}{3}j_2(\lambda\rho)\right)\delta_{ij}
-j_2(\lambda\rho)\frac{\xi_i\xi_j}{\|\boldsymbol{\xi}\|^2}\right]
\,\mathrm{d}\Phi_2(\lambda).
\end{aligned}
\]
Using the formula
\[
\frac{j_1(x)}{x}=\frac{1}{3}(j_0(x)+j_2(x)),
\]
we see that our spectral expansion is equivalent to (\ref{eq:18}) up to a constant.

To obtain the spectral representation of the field $u(\mathbf{x})$, do the following. Replace $\boldsymbol{\xi}$ with $\mathbf{x}-\mathbf{y}$ in (\ref{eq:19}), write the plane wave expansion (\ref{eq:20}) in the following form:
\begin{equation}\label{eq:32}
\begin{aligned}
e^{\mathrm{i}(\mathbf{p},\mathbf{x})}&=4\pi\sum^{\infty}_{\ell=0}
\mathrm{i}^{\ell}j_{\ell}(\lambda\rho_{\mathbf{x}})\sum^{\ell}_{m=-\ell}S^m_{\ell}
(\theta_{\mathbf{x}},\varphi_{\mathbf{x}})
S^m_{\ell}(\theta_{\mathbf{p}},\varphi_{\mathbf{p}})\\
e^{-\mathrm{i}(\mathbf{p},\mathbf{y})}&=4\pi\sum^{\infty}_{\ell'=0}
\mathrm{i}^{-\ell'}j_{\ell'}(\lambda\rho_{\mathbf{y}})\sum^{\ell'}_{m'=-\ell'}S^{m'}_{\ell'}
(\theta_{\mathbf{y}},\varphi_{\mathbf{y}})
S^{m'}_{\ell'}(\theta_{\mathbf{p}},\varphi_{\mathbf{p}}),
\end{aligned}
\end{equation}
and substitute both formulas to the modified equation (\ref{eq:19}). To simplify the result, use the following \emph{Gaunt integral} named after Gaunt \cite{Gaunt}.
\begin{equation}\label{eq:33}
\int_{S^2}S^{m_1}_{\ell_1}(\theta,\varphi)S^{m_2}_{\ell_2}(\theta,\varphi)
S^{m_3}_{\ell_3}(\theta,\varphi)\,\mathrm{d}\Omega=\sqrt{\frac{(2\ell_1+1)
(2\ell_2+1)}{4\pi(2\ell_3+1)}}g^{m_3[m_1,m_2]}_{\ell_3[\ell_1,\ell_2]}
g^{0[0,0]}_{\ell_3[\ell_1,\ell_2]}.
\end{equation}
This formula is proved in exactly the same way as in the complex case, see, for example, \cite[Proposition~3.43]{Marinucci}.

The result takes the form
\[
\begin{aligned}
R_{ij}(\mathbf{x},\mathbf{y})&=4\pi\sum^{\infty}_{\ell,\ell'=0}\sum^{\ell}_{m=-\ell}
\sum^{\ell'}_{m=-\ell'}b^{\ell'm'j}_{\ell mi,1}S^m_{\ell}
(\theta_{\mathbf{x}},\varphi_{\mathbf{x}})S^{m'}_{\ell'}
(\theta_{\mathbf{y}},\varphi_{\mathbf{y}})\\
&\quad\times\int^{\infty}_0j_{\ell}(\lambda r_{\mathbf{x}})j_{\ell'}(\lambda r_{\mathbf{y}})\,\mathrm{d}\Phi_1(\lambda)\\
&\quad+4\pi\sum^{\infty}_{\ell,\ell'=0}\sum^{\ell}_{m=-\ell}
\sum^{\ell'}_{m=-\ell'}b^{\ell'm'j}_{\ell mi,2}S^m_{\ell}
(\theta_{\mathbf{x}},\varphi_{\mathbf{x}})S^{m'}_{\ell'}
(\theta_{\mathbf{y}},\varphi_{\mathbf{y}})\\
&\quad\times\int^{\infty}_0j_{\ell}(\lambda r_{\mathbf{x}})j_{\ell'}(\lambda r_{\mathbf{y}})\,\mathrm{d}\Phi_2(\lambda).
\end{aligned}
\]
Theorem~\ref{th:1} follows from this equation and Kahrunen's theorem.

\subsection{Proof of Theorem~\ref{th:2}}

In Theorem~\ref{th:2} we have $r=2$ and $L=\mathsf{S}^2(\mathbb{R}^3)$. The coupled basis of the space $\mathsf{S}^2(\mathbb{R}^3)\otimes\mathsf{S}^2(\mathbb{R}^3)$ contains $36$ rank $4$ tensors, the tensor products of all possible pairs of the $6$ Godunov--Gordienko matrices $T^{0,1,0}_{ij}=g^{0[i,j]}_{0[1,1]}$ and $T^{2,1,n}_{ij}=g^{n[i,j]}_{2[1,1]}$, $-2\leqslant n\leqslant 2$. The uncoupled basis of the symmetric tensor product $\mathsf{S}^2(\mathsf{S}^2(\mathbb{R}^3))$ contains the 21 symmetric rank $4$ tensors shown in Table~\ref{tab:1}.

\begin{table}
\caption{The uncoupled basis of the space $\mathsf{S}^2(\mathsf{S}^2(\mathbb{R}^3))$}
\begin{tabular}{|l|l|}
\hline Tensor & Value \\
\hline
$\mathsf{T}_{ij\ell m}^{0,1,0}$ & $g^{0[i,j]}_{0[1,1]}g^{0[\ell,m]}_{0[1,1]}$ \\
$\mathsf{T}_{ij\ell m}^{0,2,0}$ & $\frac{1}{\sqrt{5}}\sum^2_{n=-2}g^{n[i,j]}_{2[1,1]}g^{n[\ell,m]}_{2[1,1]}$ \\
$\mathsf{T}_{ij\ell m}^{2,1,t}$, $-2\leqslant t\leqslant 2$ & $\frac{1}{\sqrt{6}}(\delta_{ij}g^{t[\ell,m]}_{2[1,1]}
+\delta_{\ell m}g^{t[i,j]}_{2[1,1]})$ \\
$\mathsf{T}_{ij\ell m}^{2,2,t}$, $-2\leqslant t\leqslant 2$ & $\sum^2_{n,q=-2}g^{t[n,q]}_{2[2,2]}g^{n[i,j]}_{2[1,1]}g^{q[\ell,m]}_{2[1,1]}$ \\
$\mathsf{T}_{ij\ell m}^{4,1,t}$, $-4\leqslant t\leqslant 4$ & $\sum^2_{n,q=-2}g^{t[n,q]}_{4[2,2]}g^{n[i,j]}_{2[1,1]}g^{q[\ell,m]}_{2[1,1]}$ \\
\hline
\end{tabular}
\label{tab:1}
\end{table}

By (\ref{eq:22}), the expected value of the random field $\boldsymbol{\tau}(\mathbf{x})$ is
\[
E_{ij}(\mathbf{x})=C_1T^{0,1,0}_{ij}=C\delta_{ij}.
\]

We represent the symmetric tensor $f_{ij\ell m}$ in the \emph{Voigt form} as a symmetric $6\times 6$ matrix, where Voigt indexes are numbered in the following order: $-1-1$, $00$, $11$, $01$, $-11$, $-10$. For example, $f_{-1-101}$ simplifies to $f_{14}$, and so on.

Using the values of the Godunov--Gordienko coefficients calculated in \cite{Godunov,Malyarenko}, we prove that the only non-zero elements of the symmetric matrix $f_{ij}(\mathbf{0})$ lying on and over its main diagonal are as follows:
\begin{equation}\label{eq:35}
\begin{aligned}
f_{11}(\mathbf{0})&=f_{22}(\mathbf{0})=f_{33}(\mathbf{0})=\frac{1}{3}f^{0,1,0}(0)
+\frac{2}{3\sqrt{5}}f^{0,2,0}(0),\\
f_{12}(\mathbf{0})&=f_{13}(\mathbf{0})=f_{23}(\mathbf{0})=\frac{1}{3}f^{0,1,0}(0)
-\frac{1}{3\sqrt{5}}f^{0,2,0}(0),\\
f_{44}(\mathbf{0})&=f_{55}(\mathbf{0})=f_{66}(\mathbf{0})=\frac{1}{2\sqrt{5}}f^{0,2,0}(0).
\end{aligned}
\end{equation}
It is not difficult to prove that the above matrix is nonnegative-definite with unit trace if and only if $f^{0,1,0}(0)$ and $f^{0,2,0}(0)$ are nonnegative real numbers with
\begin{equation}\label{eq:28}
f^{0,1,0}(0)+\frac{7}{2\sqrt{5}}f^{0,2,0}(0)=1.
\end{equation}

By (\ref{eq:15}), he only non-zero elements of the symmetric matrix $f_{ij}(\lambda)$ lying on and over its main diagonal are as follows:
\[
\begin{aligned}
f_{11}(\lambda)&=f_{33}(\lambda)=\frac{1}{3}f_1(\lambda)
+\frac{2}{3\sqrt{5}}f_2(\lambda)-\frac{1}{3}f_3(\lambda)
-\frac{\sqrt{2}}{3\sqrt{7}}f_4(\lambda)+\frac{3}{2\sqrt{70}}f_5(\lambda),\\
f_{12}(\lambda)&=f_{23}(\lambda)=\frac{1}{3}f_1(\lambda)
-\frac{1}{3\sqrt{5}}f_2(\lambda)+\frac{1}{6}f_3(\lambda)
-\frac{\sqrt{2}}{3\sqrt{7}}f_4(\lambda)-\frac{\sqrt{2}}{\sqrt{35}}f_5(\lambda),\\
f_{13}(\lambda)&=\frac{1}{3}f_1(\lambda)-\frac{1}{3\sqrt{5}}f_2(\lambda)
-\frac{1}{3}f_3(\lambda)+\frac{2\sqrt{2}}{3\sqrt{7}}f_4(\lambda)
+\frac{1}{2\sqrt{70}}f_5(\lambda),\\
f_{22}(\lambda)&=\frac{1}{3}f_1(\lambda)+\frac{2}{3\sqrt{5}}f_2(\lambda)
+\frac{2}{3}f_3(\lambda)+\frac{2\sqrt{2}}{3\sqrt{7}}f_4(\lambda)
+\frac{2\sqrt{2}}{\sqrt{35}}f_5(\lambda),\\
f_{44}(\lambda)&=f_{66}(\lambda)=\frac{1}{2\sqrt{5}}f_2(\lambda)
+\frac{1}{2\sqrt{14}}f_4(\lambda)-\frac{\sqrt{2}}{\sqrt{35}}f_5(\lambda),\\
f_{55}(\lambda)&=\frac{1}{2\sqrt{5}}f_2(\lambda)
-\frac{1}{\sqrt{14}}f_4(\lambda)+\frac{1}{2\sqrt{70}}f_5(\lambda).
\end{aligned}
\]
Here we introduce notation $f^{i,j,0}(\lambda)=f_{i+j}(\lambda)$. Note that $f_{13}(\lambda)=f_{11}(\lambda)-2f_{55}(\lambda)$, while $f_{12}(\lambda)$ is not a linear combination of the diagonal elements of the matrix $f_{ij}(\lambda)$. Introduce the following notation:
\begin{equation}\label{eq:34}
\begin{aligned}
u_1(\lambda)&=2f_{44}(\lambda),\qquad u_2(\lambda)=3f_{55}(\lambda),\qquad u_3(\lambda)=2(f_{11}(\lambda)-f_{55}(\lambda)),\\
u_4(\lambda)&=f_{22}(\lambda),\qquad u_5(\lambda)=f_{12}(\lambda).
\end{aligned}
\end{equation}
Direct calculations show that the matrix $f(\lambda)$ is nonnegative-definite with unit trace if and only if $u_i(\lambda)\geqslant 0$, $1\leqslant i\leqslant 4$, $u_1(\lambda)+\cdots+u_4(\lambda)=1$ and $|u_5(\lambda)|\leqslant\sqrt{u_3(\lambda)u_4(\lambda)/2}$. It follows from (\ref{eq:35}) and (\ref{eq:34}) that
\begin{equation}\label{eq:29}
u_1(0)=\frac{1}{\sqrt{5}}f_2(0),\quad u_2(0)=\frac{3}{2\sqrt{5}}f_2(0),\quad u_3(0)+u_4(0)=1-\frac{\sqrt{5}}{2}f_2(0).
\end{equation}

Define
\begin{equation}\label{eq:24}
v_1(\lambda)=\frac{u_3(\lambda)}{u_3(\lambda)+u_4(\lambda)},\qquad v_2(\lambda)=\frac{u_5(\lambda)}{u_3(\lambda)+u_4(\lambda)},
\end{equation}
and $v_1(\lambda)=1/2$, $v_2(\lambda)=0$ if the denominator is equal to $0$. We see that the set of extreme points of the set $\mathcal{C}$ contains $3$ connected components: the matrix $D^1$ with nonzero entries $D^1_{44}=D^1_{66}=1/2$, the matrix $D^2$ with nonzero entries $D^2_{11}=D^2_{33}=D^2_{55}=1/3$ and $D^2_{13}=D^2_{31}=-1/3$, and the symmetric matrices $D(\lambda)$ with nonzero entries on and over the main diagonal as follows
\[
D_{11}(\lambda)=D_{33}(\lambda)=D_{13}(\lambda)=v_1(\lambda)/2,\quad D_{22}(\lambda)=1-v_1(\lambda), \quad D_{12}(\lambda)=D_{23}(\lambda)=v_2(\lambda)
\]
lying on the ellipse
\[
u_1(\lambda)=u_2(\lambda)=0,\qquad 4(v_1(\lambda)-1/2)^2+8v^2_2(\lambda)=1.
\]
The matrix $f(\lambda)$ takes the form
\[
f(\lambda)=u_1(\lambda)D^1+u_2(\lambda)D^2+(u_3(\lambda)+u_4(\lambda))D(\lambda),
\]
where $D(\lambda)$ lies in the elliptic region $4(v_1(\lambda)-1/2)^2+8v^2_2(\lambda)\leqslant 1$.

The functions $f_i(\lambda)$ are expressed in terms of $u_i(\lambda)$ as follows:
\begin{equation}\label{eq:23}
\begin{aligned}
f_1(\lambda)&=\frac{2}{3}u_3(\lambda)+\frac{1}{3}u_4(\lambda)+\frac{4}{3}u_5(\lambda),\\
f_2(\lambda)&=\frac{2}{\sqrt{5}}u_1(\lambda)+\frac{4}{3\sqrt{5}}u_2(\lambda)
+\frac{1}{3\sqrt{5}}u_3(\lambda)+\frac{2}{3\sqrt{5}}u_4(\lambda)-\frac{4}{3\sqrt{5}}u_5(\lambda),\\
f_3(\lambda)&=-\frac{2}{3}u_3(\lambda)+\frac{2}{3}u_4(\lambda)+\frac{2}{3}u_5(\lambda),\\
f_4(\lambda)&=\frac{\sqrt{2}}{\sqrt{7}}u_1(\lambda)-\frac{4\sqrt{2}}{3\sqrt{7}}u_2(\lambda)
+\frac{\sqrt{2}}{3\sqrt{7}}u_3(\lambda)+\frac{2\sqrt{2}}{3\sqrt{7}}u_4(\lambda)
-\frac{4\sqrt{2}}{3\sqrt{7}}u_5(\lambda),\\
f_5(\lambda)&=-\frac{4\sqrt{2}}{\sqrt{35}}u_1(\lambda)+\frac{2\sqrt{2}}{3\sqrt{35}}u_2(\lambda)
+\frac{\sqrt{2}}{\sqrt{35}}u_3(\lambda)+\frac{2\sqrt{2}}{\sqrt{35}}u_4(\lambda)
-\frac{4\sqrt{2}}{\sqrt{35}}u_5(\lambda).
\end{aligned}
\end{equation}
Denote $M^{i+j}_{ij\ell m}(\mathbf{p})=M^{2i,j}_{ij\ell m}(\mathbf{p})$. By (\ref{eq:17}),
\[
f_{ij\ell m}(\mathbf{p})=M^1_{ij\ell m}(\mathbf{p})f_1(\lambda)+\cdots+M^5_{ij\ell m}(\mathbf{p})f_5(\lambda).
\]
Using (\ref{eq:24}) and (\ref{eq:23}), we obtain
\begin{equation}\label{eq:27}
\begin{aligned}
f_{ij\ell m}(\mathbf{p})&=\left[\frac{2}{\sqrt{5}}M^2_{ij\ell m}(\mathbf{p})+\frac{\sqrt{2}}{\sqrt{7}}M^4_{ij\ell m}(\mathbf{p})-\frac{4\sqrt{2}}{\sqrt{35}}M^5_{ij\ell m}(\mathbf{p})\right]u_1(\lambda)\\
&\quad+\left[\frac{4}{3\sqrt{5}}M^2_{ij\ell m}(\mathbf{p})-\frac{4\sqrt{2}}{3\sqrt{7}}M^4_{ij\ell m}(\mathbf{p})+\frac{2\sqrt{2}}{3\sqrt{35}}M^5_{ij\ell m}(\mathbf{p})\right]u_2(\lambda)\\
&\quad+\left[\frac{v_1(\lambda)+4v_2(\lambda)+1}{3}M^1_{ij\ell m}(\mathbf{p})
+\frac{-v_1(\lambda)-4v_2(\lambda)+2}{3\sqrt{5}}M^2_{ij\ell m}(\mathbf{p})\right.\\
&\quad+\frac{-4v_1(\lambda)+2v_2(\lambda)+2}{3}M^3_{ij\ell m}(\mathbf{p})
+\frac{\sqrt{2}(-v_1(\lambda)-4v_2(\lambda)+2)}{3\sqrt{7}}M^4_{ij\ell m}(\mathbf{p})\\
&\quad+\left.\frac{\sqrt{2}(-v_1(\lambda)-4v_2(\lambda)+2)}{\sqrt{35}}M^5_{ij\ell m}(\mathbf{p})\right](u_3(\lambda)+u_4(\lambda)).
\end{aligned}
\end{equation}

Substitute (\ref{eq:27}) and (\ref{eq:12}) to (\ref{eq:7}), write the result in terms of spherical harmonics (\ref{eq:16}) and use the plane wave expansion (\ref{eq:20}). We have
\begin{equation}\label{eq:30}
\begin{aligned}
R_{ij\ell m}(\boldsymbol{\xi})&=\int^{\infty}_0\left[\frac{2}{\sqrt{5}}j_0(\lambda\rho)M^2_{ij\ell m}(\boldsymbol{\xi})-\frac{\sqrt{2}}{\sqrt{7}}j_2(\lambda\rho)M^4_{ij\ell m}(\boldsymbol{\xi})\right.\\
&\quad\left.-\frac{4\sqrt{2}}{\sqrt{35}}j_4(\lambda\rho)M^5_{ij\ell m}(\boldsymbol{\xi})\right]\,\mathrm{d}\Phi_1(\lambda)\\
&\quad+\int^{\infty}_0\left[\frac{4}{3\sqrt{5}}j_0(\lambda\rho)M^2_{ij\ell m}(\boldsymbol{\xi})+\frac{4\sqrt{2}}{3\sqrt{7}}j_2(\lambda\rho)M^4_{ij\ell m}(\boldsymbol{\xi})\right.\\
&\quad\left.+\frac{2\sqrt{2}}{3\sqrt{35}}j_4(\lambda\rho)M^5_{ij\ell m}(\boldsymbol{\xi})\right]\,\mathrm{d}\Phi_2(\lambda)\\
&\quad+\int^{\infty}_0\left[\frac{v_1(\lambda)+4v_2(\lambda)+1}{3}j_0(\lambda\rho)M^1_{ij\ell m}(\boldsymbol{\xi})\right.\\
&\quad+\frac{-v_1(\lambda)-4v_2(\lambda)+2}{3\sqrt{5}}j_0(\lambda\rho)M^2_{ij\ell m}(\boldsymbol{\xi})\\
&\quad-\frac{-4v_1(\lambda)+2v_2(\lambda)+2}{3}j_2(\lambda\rho)M^3_{ij\ell m}(\boldsymbol{\xi})\\
&\quad-\frac{\sqrt{2}(-v_1(\lambda)-4v_2(\lambda)+2)}{3\sqrt{7}}j_2(\lambda\rho)M^4_{ij\ell m}(\boldsymbol{\xi})\\
&\quad+\left.\frac{\sqrt{2}(-v_1(\lambda)-4v_2(\lambda)+2)}{\sqrt{35}}j_4(\lambda\rho)M^5_{ij\ell m}(\boldsymbol{\xi})\right]\,\mathrm{d}\Phi_3(\lambda),
\end{aligned}
\end{equation}
where we introduced notation $\mathrm{d}\Phi_j(\lambda)=u_j(\lambda)\,\mathrm{d}\mu(\lambda)$, $j=1$, $2$, and $\mathrm{d}\Phi_3(\lambda)=(u_3(\lambda)+u_4(\lambda))\,\mathrm{d}\mu(\lambda)$. It follows from (\ref{eq:28}) that $0\leqslant f_2(0)\leqslant\frac{2\sqrt{5}}{7}$. Then, by (\ref{eq:29}),
\[
\frac{2}{7}\leqslant u_3(0)+u_4(0)\leqslant 1.
\]
It follows that the atom $\Phi_3(\{0\})$ occupies at least $2/7$ of the sum of all three atoms, while the rest is divided between $\Phi_1(\{0\})$ and $\Phi_2(\{0\})$ in the proportion $1:\frac{3}{2}$.

In \cite{Malyarenko} we proved that $M^n_{ij\ell m}(\mathbf{p})$ are expressed in terms of $L^n_{ij\ell m}(\mathbf{p})$ as follows.
\begin{equation}\label{eq:25}
\begin{aligned}
M^1_{ij\ell m}(\mathbf{p})&=\frac{1}{3}L^1_{ij\ell
m}(\mathbf{p}),\\
M^2_{ij\ell
m}(\mathbf{p})&=-\frac{1}{3\sqrt{5}}L^1_{ij\ell m}(\mathbf{p})
+\frac{1}{2\sqrt{5}}L^2_{ij\ell m}(\mathbf{p}),\\
M^3_{ij\ell
m}(\mathbf{p})&=-\frac{1}{3}L^1_{ij\ell m}(\mathbf{p})
+\frac{1}{2}L^4_{ij\ell m}(\mathbf{p}),\\
M^4_{ij\ell
m}(\mathbf{p})&=\frac{2\sqrt{2}}{3\sqrt{7}}L^1_{ij\ell m}(\mathbf{p})
-\frac{1}{\sqrt{14}}L^2_{ij\ell m}(\mathbf{p})
+\frac{3}{2\sqrt{14}}L^3_{ij\ell m}(\mathbf{p})-
\frac{\sqrt{2}}{\sqrt{7}}L^4_{ij\ell m}(\mathbf{p}),\\
M^5_{ij\ell
m}(\mathbf{p})&=\frac{1}{2\sqrt{70}}L^1_{ij\ell m}(\mathbf{p})
-\frac{13}{2\sqrt{70}}L^2_{ij\ell m}(\mathbf{p})
-\frac{\sqrt{5}}{2\sqrt{14}}L^3_{ij\ell m}(\mathbf{p})
-\frac{\sqrt{5}}{2\sqrt{14}}L^4_{ij\ell m}(\mathbf{p})\\
&\quad+\frac{\sqrt{35}}{2\sqrt{2}}L^5_{ij\ell m}(\mathbf{p}).\nonumber
\end{aligned}
\end{equation}
The second and fourth equations were proved by brutal force, using the values of matrix entries and Godunov--Gordienko coefficients calculated in \cite{Godunov,Malyarenko}. Here is the algebraic proof.

It follows from the definition of the Godunov--Gordienko coefficients that
\begin{equation}\label{eq:26}
D^{\ell_1}_{m_1n_1}(k)D^{\ell_2}_{m_2n_2}(k)=\sum^{\ell_1+\ell_2}_{\ell=|\ell_1-\ell_2|}
\sum^{\ell}_{q_1,q_2=-\ell}g^{q_1[m_1,m_2]}_{\ell[\ell_1,\ell_2]}D^{\ell}_{q_1q_2}(k)
g^{q_2[n_1,n_2]}_{\ell[\ell_1,\ell_2]}.
\end{equation}
Put $k=I$, $m_1=i$, $m_2=j$, $n_1=\ell$, $n_2=m$, and $\ell_1=\ell_2=1$. We obtain
\[
\delta_{i\ell}\delta_{jm}=\frac{1}{3}\delta_{ij}\delta_{\ell m}+\sum^1_{n=-1}g^{n[i,j]}_{1[1,1]}g^{n[\ell,m]}_{1[1,1]}+\sqrt{5}M^2_{ij\ell m}(\mathbf{p}).
\]
Interchange $\ell$ and $m$ and use the fact that $g^{n[\ell,m]}_{1[1,1]}$ is a skew-symmetric matrix. Then
\[
\delta_{im}\delta_{j\ell}=\frac{1}{3}\delta_{ij}\delta_{\ell m}-\sum^1_{n=-1}g^{n[i,j]}_{1[1,1]}g^{n[\ell,m]}_{1[1,1]}+\sqrt{5}M^2_{ij\ell m}(\mathbf{p}).
\]
Adding two last displays yields
\[
L^2_{ij\ell m}(\mathbf{p})=\frac{2}{3}L^1_{ij\ell m}(\mathbf{p})+2\sqrt{5}M^2_{ij\ell m}(\mathbf{p}),
\]
which is equivalent to the second equation in (\ref{eq:25}).

It is proved in \cite{Malyarenko} that
\[
\frac{p_ip_j}{\|\mathbf{p}\|^2}=\frac{1}{3}\delta_{ij}+\sqrt{2/3}\sum^2_{n=-2}
g^{n[i,j]}_{2[1,1]}D^2_{n0}(\mathbf{p}).
\]
Rewrite this equation as
\[
\frac{p_{\ell}p_m}{\|\mathbf{p}\|^2}=\frac{1}{3}\delta_{\ell m}+\sqrt{2/3}\sum^2_{q=-2}
g^{q[\ell,m]}_{2[1,1]}D^2_{q0}(\mathbf{p})
\]
and multiply both equations. We obtain
\[
L^5_{ij\ell m}(\mathbf{p})=\frac{1}{9}L^1_{ij\ell m}(\mathbf{p})+\frac{3}{2}
M^3_{ij\ell m}(\mathbf{p})+\frac{2}{3}\sum^2_{n,q=-2}g^{n[i,j]}_{2[1,1]}g^{q[\ell,m]}_{2[1,1]}
D^2_{n0}(\mathbf{p})D^2_{q0}(\mathbf{p}).
\]
By (\ref{eq:26}),
\[
D^2_{n0}(\mathbf{p})D^2_{q0}(\mathbf{p})=g^{0[n,q]}_{0[2,2]}g^{0[0,0]}_{0[2,2]}
+g^{0[0,0]}_{2[2,2]}\sum^2_{s=-2}g^{s[n,q]}_{2[2,2]}D^2_{s0}(\mathbf{p})
+g^{0[0,0]}_{4[2,2]}\sum^4_{s=-4}g^{s[n,q]}_{4[2,2]}D^4_{s0}(\mathbf{p}).
\]
Using the values $g^{0[n,q]}_{0[2,2]}=\sqrt{1/5}\delta_{nq}$, $g^{0[0,0]}_{2[2,2]}=\sqrt{2/7}$, and $g^{0[0,0]}_{4[2,2]}=\frac{3\sqrt{2}}{\sqrt{35}}$ calculated in \cite{Malyarenko}, after simple algebraic calculations we obtain the fourth equation in (\ref{eq:25}).

Apply (\ref{eq:16}) and (\ref{eq:25}) to (\ref{eq:30}). We obtain (\ref{eq:31}).

\subsection{Proof of Theorem~\ref{th:3}}

Substitute (\ref{eq:27}) and (\ref{eq:12}) to (\ref{eq:7}), write the result in terms of spherical harmonics (\ref{eq:16}), replace $\boldsymbol{\xi}$ with $\mathbf{x}-\mathbf{y}$ and use the plane wave expansion in the form (\ref{eq:32}). To simplify the result, use the Gaunt integral (\ref{eq:33}).

\section{Concluding remarks}\label{sec:4}

Methods of our paper work equally good in the case of $r=0$. The convex compact set $\mathcal{C}$ is a one-point set, and one can deduce the classical results by Schoenberg \cite{Schoenberg}
\[
R(\boldsymbol{\xi})=\int^{\infty}_0\frac{\sin(\lambda\|\boldsymbol{\xi}\|)}
{\lambda\|\boldsymbol{\xi}\|}\,\mathrm{d}\Phi(\lambda)
\]
and Yadrenko \cite{Yadrenko}
\[
T(\rho,\theta,\varphi)=C+2\sqrt{\pi}\sum^{\infty}_{\ell=0}\sum^{\ell}_{m=-\ell}
\int^{\infty}_0j_{\ell}(\lambda\rho)\,\mathrm{d}Z^m_{\ell}(\lambda)S^m_{\ell}
(\theta,\varphi).
\]

Consider a particular case of Theorems~\ref{th:2} and \ref{th:3} when $u_5(\lambda)=0$. It means that the random fields $\tau_{00}(\mathbf{x})$ and $(\tau_{-1-1}(\mathbf{x}),\tau_{11}(\mathbf{x}))^{\top}$ are uncorrelated. In this case the set $\mathcal{C}$ becomes a tetrahedron with four extreme points: $D^1$, $D^2$, $D^3$, and $D^4$. The matrix $D^3$ is equal to $D(\lambda)$ when
\begin{equation}\label{eq:36}
v_1(\lambda)=1,\qquad v_2(\lambda)=0,
\end{equation}
while the matrix $D^4$ is equal to $D(\lambda)$ when
\begin{equation}\label{eq:37}
v_1(\lambda)=v_2(\lambda)=0.
\end{equation}
The spectral expansion of Theorem~\ref{th:2} takes the form
\[
B_{ij\ell m}(\boldsymbol{\xi})=\sum^4_{n=1}\int^{\infty}_0\sum^5_{q=1}
N_{nq}(\lambda,\rho)L^q_{ij\ell m}(\boldsymbol{\xi})\,\mathrm{d}\Phi_n(\lambda),
\]
where $N_3q(\lambda,\rho)$ (resp. $N_4q(\lambda,\rho)$) can be calculated by substituting (\ref{eq:36}) (resp. (\ref{eq:37})) to the last five elements of the third column of Table~\ref{tab:2}. If
\[
\sum^4_{n=1}\Phi_n(\{0\})=\Phi_0>0,
\]
then
\begin{eqnarray*}
\Phi_1(\{0\})&=&\frac{\Phi_0f_2(0)}{\sqrt{5}},\qquad
\Phi_2(\{0\})=\frac{3\Phi_0f_2(0)}{2\sqrt{5}},\\
\Phi_3(\{0\})&=&\Phi_0\left(\frac{1}{3}-\frac{f_2(0)}{2\sqrt{5}}\right),\qquad
\Phi_4(\{0\})=\Phi_0\left(\frac{2}{3}-\frac{2f_2(0)}{\sqrt{5}}\right),
\end{eqnarray*}
with $0\leqslant f_2(0)\leqslant 2\sqrt{5}/7$. The spectral expansion of Theorem~\ref{th:3} takes the form
\[
\tau_{ij}(r,\theta,\varphi)=C\delta_{ij}+2\sqrt{\pi}\sum^4_{n=1}
\sum^{\infty}_{u=0}\sum^u_{w=-u}\int^{\infty}_0j_u(\lambda r)\,\mathrm{d}Z^{n'}_{uwij}(\lambda)S^w_u(\theta,\varphi),
\]
where the measures $\mathrm{d}Z^{n'}_{uwij}(\lambda)$ are determined by (\ref{eq:38}). In (\ref{eq:38}), $L^n$ are infinite lower triangular matrices from Cholesky factorisation of nonnegative-definite matrices $b^{u'w'\ell m}_{uwij,n}$. The matrix $b^{u'w'\ell m}_{uwij,3}$ (resp. $b^{u'w'\ell m}_{uwij,4}$) can be calculated by substituting (\ref{eq:36}) (resp. (\ref{eq:37})) to the formula that determines $b^{u'w'\ell m}_{uwij,3}(\lambda)$.

We conjecture that in the general case the set of extreme points of the convex compact set $\mathcal{C}$ has finitely many, say $N$, connected components. Each of the components is either an one-point set or an ellipsoid. To each connected component $\mathcal{D}_i$ we associate a pair $(\Phi_i,\mathbf{v}_i)$, where $\Phi_i$ is a finite measure on $[0,\infty)$, and $\mathbf{v}_i$ is a $\Phi_i$-equivalence class of measurable functions on $[0,\infty)$ with values in the closed convex hull of the set $\mathcal{D}_i$ (constant if $\mathcal{D}_i$ is a one-point set). The number of integrals in the spectral representation is equal to $N$, and the $i$th integral is taken with respect to the measure $\Phi_i$.

\end{document}